\documentclass[12pt,centertags,oneside]{amsart}
\usepackage{amsmath,amstext,amsthm,amscd,typearea,hyperref}
\usepackage{amssymb}
\usepackage{a4wide}
\usepackage[mathscr]{eucal}
\usepackage{mathrsfs}
\usepackage{typearea}
\usepackage{charter}
\usepackage{pdfsync}
\usepackage[a4paper,width=16.2cm,top=3cm,bottom=3cm]{geometry}

\numberwithin{equation}{section}

\newtheorem{theorem}{Theorem}[section]
\newtheorem{definition}[theorem]{Definition}
\newtheorem{proposition}[theorem]{Proposition}
\newtheorem{corollary}[theorem]{Corollary}
\newtheorem{lemma}[theorem]{Lemma}
\newtheorem{remark}[theorem]{Remark}

\newtheorem{example}[theorem]{Example}

\newcommand{\cali}[1]{\mathscr{#1}}

\newcommand{\supp}{{\rm supp}}

\newcommand{\dist}{\mathop{\mathrm{dist}}\nolimits}

\newcommand{\vol}{\mathop{\mathrm{vol}}}
\newcommand{\sgn}{\mathop{\mathrm{sgn}}}
\newcommand{\eq}{\mathrm{eq} }

\newcommand{\ddc}{dd^c}

\def\mE{\mathcal{E}}

\def\mL{\mathcal{L}}

\newcommand{\PSH}{{\rm PSH}}

\newcommand{\Sym}{{\rm Sym}}

\newcommand{\Bc}{\cali{B}}
\newcommand{\Cc}{\cali{C}}

\renewcommand{\Mc}{\cali{M}}
\newcommand{\Oc}{\cali{O}}
\newcommand{\Pc}{\cali{P}}

\newcommand{\C}{\mathbb{C}}
\newcommand{\D}{\mathbb{D}}

\newcommand{\R}{\mathbb{R}}
\renewcommand\P{\mathbb{P}}

\renewcommand{\S}{\mathbb{S}}

\newcommand{\x}{{\bfit{x}}}
\newcommand{\y}{{\bfit{y}}}

\newcommand{\bfit}[1]{\textbf{\textit{#1}}}

\title[Large deviations principle for  some Beta ensembles]{ 
Large deviations principle for some Beta ensembles}

\author{Tien-Cuong Dinh}
\address{Department of Mathematics, National University 
of Singapore, 10 Lower Kent Ridge Road, Singapore 119076.}
\email{matdtc@nus.edu.sg}
 
\author{Vi{\^e}t-Anh Nguy{\^e}n}
\address{Math{\'e}matique-B{\^a}timent 425, UMR 8628, 
Universit{\'e} Paris-Sud, 91405 Orsay, France. 
}
\email{VietAnh.Nguyen@math.u-psud.fr}
\date{March 11, 2016}
\begin{document}

\begin{abstract} Let $L$ be a positive line bundle over a projective complex manifold $X$, $L^p$ its tensor power of order $p$,  $H^0(X,L^p)$ the space of holomorphic sections of $L^p$ and $N_p$ the complex dimension of $H^0(X,L^p)$. 
The determinant of a basis of $H^0(X,L^p)$, together with some given probability measure on a  weighted compact set in $X$, induces naturally 
a $\beta$-ensemble, i.e.,  a  random $N_p$-point  process  on the compact set. 
Physically, this general setting corresponds to a gas of free fermions on $X$ and may admit some random matrix models.
The empirical measures, associated  with such $\beta$-ensembles, converge almost surely to an equilibrium measure  when $p$ goes to infinity.  
We  establish  a  large  deviations principle  
(LDP)  with  an effective speed of convergence
 for   these  empirical  measures. 
Our study covers the case of some $\beta$-ensembles on a compact subset of the unit sphere $\S^n\subset \R^{n+1}$ or of the Euclidean space $\R^n$. 
\end{abstract}

\maketitle

\medskip

\noindent
{\bf Classification AMS 2010}: 32U15 (32L05, 60F10).

\medskip




\section{Introduction} \label{introduction}

Let $K$ be a metric space and  $N$  a positive integer. If $\bfit{x}=(x_1,\ldots, x_{N})$ is a point in the  $N$-fold product  $K^N,$ then the associated {\it empirical measure} is the probability measure 
$$\mu^{\bfit{x}}:={1\over N} \sum_{k=1}^N \delta_{x_k}$$
which is equidistributed on $x_1,\ldots, x_N$. 
Here, $\delta_x$ denotes the Dirac mass at $x$. 
Any  probability measure $\nu$ on $K^N$ induces a random $N$-point process on $K$  and $\nu$ is  the law of this random process.

Let $\{N_p\}_{p\geq 1}$ be a sequence of positive integers such that
 $N_p\to\infty$    as $p\to\infty$  
 and $\{\nu_p\}_{p\geq 1}$ a sequence of probability measures on $K^{N_p}$. In many problems from mathematics or mathematical physic, 
a central question is to study the eventual convergence of the sequence $\mu^{\x^{(p)}}$ to an equilibrium measure, where  $\x^{(p)}$ is the random $N_p$-point process on $K$  described  by the  law $\nu_p.$  A significantly interesting setting considered in literature is the case of $\beta$-ensembles on a compact subset of the unit sphere $\S^n$ in $\R^{n+1}$ or a compact subset of $\R^n$. We will obtain in this paper a LDP for such $\beta$-ensembles with an explicit rate of convergence.     
Our approach uses however techniques from complex analysis, and therefore we will first describe  the general setting, which, 
physically,  corresponds to a gas of free fermions and may admit some random matrix models. 
The reader will find in the paper of Berman \cite{Berman} a detailed exposition and a list of references. 
The case of $\beta$-ensembles on the unit sphere or on the real Euclidean space, mentioned above,  will be obtained as a  corollary, see Examples \ref{e:space} and \ref{e:sphere} below. 

Let $X$ be a compact K\"ahler manifold of dimension $n$. 
Let $L\to X$ be a positive line bundle endowed with a given smooth Hermitian metric $h_0$. We assume that the metric $h_0$ is positively curved, that is, the Chern form $\omega_0$ associated with $h_0$ is a K\"ahler form on $X$. 
For simplicity, we will use the Riemannian metric on $X$ induced by $\omega_0$.
The space of holomorphic sections of $L^p:=L\otimes \cdots\otimes L$ ($p$ times) is denoted by $H^0(X,L^p).$  
Since $L$ is  ample,  
by  Kodaira-Serre vanishing and Riemann-Roch-Hirzebruch theorems (see \cite[Thm 1.5.6 and 1.4.6]{MaMarinescu07}),
we have
\begin{equation}\label{e:Demailly}
N_p:=\dim H^0(X,L^p) =  {p^n\over n!} \|\omega_0^n\| +O(p^{n-1}).
\end{equation}
Here, $ \|\omega_0^n\|$ denotes the mass of the volume form $\omega_0^n$. It depends only on the Chern class of $L$. 

If $L_1, L_2$ are line bundles over complex manifolds $X_1$ and  $X_2$ respectively, we denote by $L_1 \boxtimes L_2$
the line bundle over the product manifold $X_1 \times X_2$ defined as
 $L_1 \boxtimes L_2 := \pi_1^* (L_1)\otimes \pi^*_2 (L_2),$ where $\pi_1,  \pi_2$ are the natural projections from $X_1 \times X_2$ to its factors. 
 If $L_1$ and $L_2$ are endowed with some Hermitian metrics, then 
 $L_1 \boxtimes L_2$ carries also a metric induced by those on $L_1$ and $L_2$.

Let $S_p=(s_1,\ldots,s_{N_p})$ be a basis of $H^0(X,L^p)$. We define the section $\det S_p$ of the line bundle $(L^p)^{\boxtimes N_p}:=
L^p\boxtimes\cdots\boxtimes L^p$  ($N_p$ times) over $X^{N_p}$ by the  identity 
$$\det  S_p(\x^{(p)}) := \sum_{\sigma\in  \Sym _{N_p}} \sgn(\sigma) 
\bigotimes_{i=1}^{N_p} s_i(x_{\sigma(i)} ) \qquad \text{for} \qquad \x^{(p)}=(x_1,\ldots, x_{N_p})\in X^{N_p}, $$
where $\Sym_{N_p}$ denotes the permutation group  of $\{1,\ldots,N_p\}.$
Note that when we change the basis $S_p$, this section only changes by a non-zero multiplicative constant.

Let  $K$ be a 
compact set  in $X$ and $\phi$   a continuous real-valued function on $K$. We say that the pair $(K,\phi)$ is a {\it weighted compact set}. Let  
 $\mu$ be  a probability measure on $K.$
 
\begin{definition}\label{D:beta} \rm
Let $\beta > 0$ be a constant.   {\it A $\beta$-ensemble}  associated with the line bundle $L^p$, the weighted compact set $(K,\phi)$ and the probability measure $\mu$, is the random $N_p$-point  process on
$K$ whose joint distribution is given by
\begin{equation}\label{e:beta}
\nu^\beta_p:=c_{p,\beta} \|\det S_p(\x^{(p)})\|^\beta e^{-\beta p(\phi(x_1)+\cdots+\phi(x_{N_p}))} d\mu(x_1) 
 \otimes\cdots\otimes
d\mu(x_{N_p}),
 \end{equation}
where $c_{p,\beta}$ is  the normalizing constant so that  $\nu^\beta_p$  is a probability measure on $K^{N_p}$.
\end{definition}

Observe  that the constant $c_{p,\beta}$ depends also on $L^p,K,\phi,\mu$, but the above random point process, i.e., the measure $\nu^\beta_p$,  is  independent of the  choice of the basis $S_p$ of $H^0(X,L^p)$.
We will study these $\beta$-ensembles when $p$ goes to infinity. We need some assumptions on the regularity of $K,\phi$ and $\mu$. 
Under such conditions, we will see later that the sequence $\mu^{\x^{(p)}}$ converges almost surely to a limit $\mu_\eq(K,\phi)$ which is called the equilibrium measure of the weighted compact set $(K,\phi)$.

 Recently,  Berman \cite{Berman} obtained  a LDP in the  spirit of Donsker and Varadhan 
 \cite{DemboZeitouni} using some functionals on the space of measures. 
 In the case where $K=X$, $\phi=0$, and $\mu$ is the Lebesgue measure on $X$, 
  Carroll,  Marzo,  Massaneda and   Ortega-Cerd\`a obtained precise and optimal estimates on the expectation 
  of the Kantorovich-Wasserstein distance between $\mu^{\x^{(p)}}$ and $\mu_\eq(K,\phi)$ when $p\to\infty$  \cite{Carroll}.
  An advantage of the latter work is that Kantorovich-Wasserstein distance gives us a very explicit information about the convergence of $\mu^{\x^{(p)}}$ to $\mu_\eq(K,\phi)$.
  Our aim is to establish a LDP with precise estimations in a quite general setting and in the sprit of  the work by Carroll,  Marzo,  Massaneda and   Ortega-Cerd\`a. In order to state the main result, we need to introduce some more notions.

Let  $\Mc(X)$ denote the space of all
(Borel) probability measures on $X$. 
For $\gamma>0$, define
 the distance $\dist_\gamma$ between two measures 
   $\mu$ and $\mu'$ in $\Mc(X)$ by 
$$\dist_\gamma(\mu,\mu'):=\sup_{\| v\|_{\Cc^\gamma}\leq 1} \big |\langle \mu-\mu', v\rangle\big|,$$
where $v$ is a test smooth real-valued function. 
This distance induces the weak topology on $\Mc(X)$. 
By interpolation between Banach spaces 
(see \cite{DinhSibony4, Triebel}), 
for $0<\gamma\leq \gamma'$, there exists a constant $c>0$ such that
\begin{align}\label{eq:n1.23}
\dist_{\gamma'}\leq \dist_\gamma
\leq c[\dist_{\gamma'}]^{\gamma/\gamma'}.
\end{align}
Note that $\dist_1$ is equivalent to the classical Kantorovich-Wasserstein distance.

In Section \ref{S:regularity} below,  we  will single  out
a very large class of compact sets $K$ which enjoy the so-called $(\Cc^\alpha,\Cc^{\alpha'})$-regularity.  We will also introduce 
the notion of $\delta$-Bernstein-Markov measures which enjoy
a  quantified version of the  Bernstein-Markov property. Here, $\delta$ is a constant such that $0<\delta <1.$ 
Having in hand  these natural notions, we are in the position to state  the main result of the paper.

\begin{theorem}\label{T:abstract}
Let $X$ be a complex projective manifold of dimension $n$. Let $L$ be a positive line bundle over $X$ endowed with a smooth positively curved Hermitian metric $h_0$.
Let   $\beta>0$  and $0<\gamma\leq 2$ be constants. 
 Let $K$ be  a  $(\Cc^\alpha,\Cc^{\alpha'})$-regular compact subset of $X$ and   
  $\phi$ a $\Cc^\alpha$ real-valued function on $K$ for some constants $0<\alpha\leq 2$ and $ 0<\alpha'\leq 1$. Let  $\mu$ be a probability measure on $K$ 
 which  is   $\delta$-Bernstein-Markov with respect to  $(K, \phi)$  for some $0<\delta<1.$
Then, for every $\lambda>0,$ there are $c>0$  and Borel sets $E_{p}\subset K^{N_p}$ such that 
\begin{enumerate}
\item[(a)] $\nu^\beta_p(E_{p})\leq e^{-\lambda p^{n+1-\delta} }$;
\item[(b)] if  $\mu_p^{\x}$ denotes the empirical measure associated with $\x\in K^{N_p} \setminus E_{p},$  
then   
$$\dist_\gamma(\mu^{\x}, \mu_\eq(K,\phi))\leq c q^{\gamma}.$$
\end{enumerate}
Here, 
$q:=p^{-\delta/4}$ if $\delta/4<\alpha''$, 
$q:=p^{-\alpha''}(\log p)^{3\alpha''}$ if $\delta/4\geq \alpha''$, and  $\alpha'':=\alpha'/(24+12\alpha')$.
 \end{theorem}
 
If a sequence of points $\x^{(p)} \in K^{N_p}$ satisfies $\x^{(p)}\not\in E_p$ for $p$ large enough, then we deduce from the last theorem that 
$\mu^{\x^{(p)}}\to \mu_\eq(K,\phi)$ when $p$ goes to infinity. Therefore,
$\mu^{\x^{(p)}}$ converge almost surely to $\mu_\eq(K,\phi)$ when $p$ goes to infinity. 
More precisely,  the infinite product  $\nu^\beta:=\nu_1^\beta\times \nu_2^\beta\times  \cdots$ is a probability measure on the space of all sequences $(\x^{(p)})_{p=1}^\infty$. With respect to this measure,
the convergence $\mu^{\x^{(p)}}\to \mu_\eq(K,\phi)$ holds for almost every sequence $(\x^{(p)})_{p=1}^\infty$.

The estimate on the size of $E_p$ is a version of LDP. 
Our result also implies that 
 \begin{equation}\label{e:expectation}
\int_{X^{N_p}} \dist_\gamma\big(\mu^\x, \mu_\eq(K,\phi)\big) d\nu^\beta_p(\x)=O( q^{\gamma}).
\end{equation}
This distance expectation estimate is similar to the one obtained by
 Carroll,  Marzo,  Massaneda and   Ortega-Cerd\`a in \cite{Carroll} that we mentioned above. 
These authors proved for $K=X$, $\phi=0$ and $\mu$ the normalized  Lebesgue measure on $X$ that there is a constant $c>0$ satisfying
\begin{align}\label{eq:Carroll}
 c^{-1} p^{1/2}\leq  \int_{X^{N_p}} \dist_1\big(\mu^\x, \mu_\eq(K,\phi)\big) d\nu^\beta_p(\x)\leq  c p^{1/2}
\end{align}
for all $p.$

In order to get more concrete applications of our main result, we need  the following natural class of positive Borel measures.

\begin{definition}\label{D:mass-density} \rm
We say that  a  positive measure $\mu$ on $X$ satisfies {\it the mass-density condition} with respect to
a compact $W\subset X$ if there are two constants $c>0$ and $\rho>0$ such that
$$
\mu(B(x,r)  )\geq  cr^\rho \qquad \text{ for }\quad x\in W\quad \text{and}\quad 0<r<1.
$$
Here, $B(x,r)$ denotes the ball in $(X,\omega_0)$ of radius $r$ and centered at the point $x$.
\end{definition}

Assume now that $K$ is a smooth real manifold in $X$ with piecewise smooth boundary such that the tangent space of $K$ at each point is not  contained in a complex hyperplane of the tangent space of $X$ at that point. It was shown in \cite{DinhMaNguyen, Vu},  for $0<\alpha<1$,  that $K$ is $(\Cc^\alpha,\Cc^{\alpha/2})$-regular and is  $(\Cc^\alpha,\Cc^{\alpha})$-regular  when its boundary is smooth, see Theorem \ref{t:DMNV} below. In this case, if $\mu$ is a probability measure on $K$ satisfying the above mass-density condition for $W=K$, we will show in Corollary \ref{C:Bernstein-Markov} below that it satisfies the $\delta$-Bernstein-Markov property required in Theorem \ref{T:abstract}. Therefore, the following result is a direct consequence of that theorem.

 \begin{corollary}\label{C:manifold}
 Let $X,L,h_0, \beta,\gamma$ be as in Theorem \ref{T:abstract}. 
Let $K$ be a smooth real manifold in $X$ with piecewise smooth boundary such that the tangent space of $K$ at each point is not  contained in a complex hyperplane of the tangent space of $X$ at that point. 
Let $\mu$ be  a  probability measure on $K$ satisfying the mass-density condition  with respect to $K.$
 Let   $\phi$ be a $\Cc^\alpha$ real-valued function on $K$ with $0<\alpha<1.$ 
 Then, for every $0<\delta<1,$ the conclusion of Theorem \ref{T:abstract} holds for 
 $\alpha'':=\alpha/(48+24\alpha)$. Moreover, if the boundary of $K$ is smooth, the the same statement holds for  $\alpha'':=\alpha/(24+12\alpha)$. 
 \end{corollary}
 
Of course,  Corollary  \ref{C:manifold} holds when $\mu$ is given by the normalized volume form on $K$. 
It is  worthy  noting that  the  assumption on the  mass-density condition  of  the measure $\mu$ in this result  can be  weakened.   
In fact,  we only need that $\mu$ satisfies the mass-density condition 
on a  subset $W\subset K$  which satisfies a maximum principle, see  Corollary \ref{C:Bernstein-Markov} below.

\begin{example} \rm \label{e:space}
Let $K$ be the closure of an open set with piecewise smooth boundary in $\R^{n}$. Let $\phi$ be a $\Cc^\alpha$ real-valued function on $K$ and $\mu$ a probability measure on $K$ which satisfies the mass-density condition with respect to $K$. It is already interesting to consider the case where $\mu$  is the normalization of the restriction to $K$ of the Lebesgue measure on $\R^n$.
Denote by $\Pc_p$ the set of real polynomials of degree at most $p$ and $N_p$ the dimension of $\Pc_p$.
Choose a basis $(P_1,\ldots,P_{N_p})$ of $\Pc_p$. 
Define the probability measure $\nu_p^\beta$ at a point $\x=(x_1,\ldots, x_{N_p})$ on $K^{N_p}$ by 
$$ c_{p,\beta} |\det(P_i(x_j))|^\beta e^{-\beta p(\phi(x_1)+\cdots+\phi(x_{N_p}))}e^{-{1\over 2}{\beta p} (\log(1+\|x_1\|^2) +\cdots + \log(1+\|x_{N_p}\|^2) ) }\mu(x_1)\otimes \cdots\otimes \mu(x_{N_p}),$$
where $c_{p,\beta}$ is a normalizing constant so that $\nu_p^\beta$ is a probability measure.  Here, $\det(\cdot)$ denotes the standard determinant of a square matrix. Then the conclusion of Theorem \ref{T:abstract} holds for $\alpha'':=\alpha/(48+24\alpha)$. If the boundary of $K$ is smooth, we can take $\alpha'':=\alpha/(24+12\alpha)$. The equilibrium measure $\mu_\eq(K,\phi)$ is a probability measure supported by $K$. Its definition is given in Section \ref{S:regularity}.

In order to obtain this result as a consequence of Theorem \ref{T:abstract} and Corollary \ref{C:manifold}, consider $\R^{n}$ as the real part of $\C^{n}$ and $\C^{n}$ as a Zariski open set of the projective space $\P^{n}$. Denote by $[z_0:\cdots:z_{n}]$ the homogeneous coordinates of $\P^{n}$. We identify $\C^{n}$ with the open set $\{z_0=1\}$. Define  $X:=\P^n$. We can identify, in the natural way, the polynomials of degree $\leq p$ on $\R^{n}$ with holomorphic sections of $L^p$ with $L=\Oc(1)$ the tautological line bundle of $X=\P^{n}$. We consider the standard Hermitian metrics on these line bundles.  So $\{P_1,\ldots, P_{N_p}\}$ is identified to a basis of $H^0(X,L^p)$. 
If a section $s$ in $H^0(X,L^p)$ is identified to a polynomial $P$ then 
$$\|s(z)\|=|P(z)| e^{-{1\over 2} p\log (1+\|z\|^2)} \quad \text{for} \quad z\in\C^n.$$
So the factor involving $\log(1+\|x_i\|^2)$ in the definition of $\nu^\beta_p$ is due to the standard Hermitian metric of $L^p$. We can now apply Theorem \ref{T:abstract} and Corollary \ref{C:manifold} and get the LDP in this case. An interesting particular situation is the case where the weight $\phi$ is equal to $-{1\over 2}\log(1+\|\cdot\|^2)$. 
\end{example}

\begin{example} \rm \label{e:sphere}
Let $K$ be the closure of an open set with piecewise smooth boundary in the unit sphere $\S^n$ of $\R^{n+1}$. Let $\phi$ be a $\Cc^\alpha$ real-valued function on $K$ and $\mu$ a probability measure on $K$ which satisfies the mass-density condition with respect to $K$. It is already interesting to consider the case where $\phi=0$ and $\mu$ is  the normalization of the restriction to $K$ of the Haar measure on $\S^n$.
Consider the functions which are restrictions of (real) polynomials on $\R^{n+1}$ to $\S^n$.
Denote by $\Pc_p$ the set of these functions obtained by using polynomials of degree at most $p$ and $N_p$ the dimension of $\Pc_p$.
Note that $\Pc_p$ is isomorphic to the quotient of the space of polynomials of degree $\leq p$ by the subspace of polynomials divisible by
$x_1^2+\cdots+x_{n+1}^2-1$, where $(x_1,\ldots, x_{n+1})$ is the standard coordinate system of $\R^{n+1}$. 

Choose a basis $(P_1,\ldots,P_{N_p})$ of $\Pc_p$. 
Define the probability measure $\nu_p^\beta$ on $K^{N_p}$ by 
$$\nu_p^\beta (\x):= c_{p,\beta} |\det(P_i(x_j))|^\beta e^{-\beta p(\phi(x_1)+\cdots+\phi(x_{N_p}))}\mu(x_1)\otimes \cdots\otimes \mu(x_{N_p}),$$
where $\x=(x_1,\ldots, x_{N_p})$ is a point in $K^{N_p}$ and $c_{p,\beta}$ is a normalizing constant so that $\nu_p^\beta$ is a probability measure.
Then the conclusion of Theorem \ref{T:abstract} holds for $\alpha'':=\alpha/(48+24\alpha)$. If the boundary of $K$ is smooth, we can take $\alpha'':=\alpha/(24+12\alpha)$. The measure $\mu_\eq(K,\phi)$ is supported by $K$. In the case where $K=\S^n$ and $\phi=0$, by symmetry, this measure coincides with the Haar measure on $\S^n$.

In order to obtain this result as a consequence of Theorem \ref{T:abstract} and Corollary \ref{C:manifold}, we need to complexify $\S^n$. 
Consider $\R^{n+1}$ as the real part of $\C^{n+1}$ and $\C^{n+1}$ as a Zariski open set of the projective space $\P^{n+1}$. Denote by $[z_0:\cdots:z_{n+1}]$ the homogeneous coordinates of $\P^{n+1}$. We identify $\C^{n+1}$ with the open set $\{z_0=1\}$. The sphere $\S^n$ is then the intersection of $\R^{n+1}$ with the complex hypersurface $z_1^2+\cdots+z_{n+1}^2=z_0^2$ in $\P^{n+1}$. Denote by $X$ this hypersurface. We can identify, in the natural way, the polynomials of degree $\leq p$ on $\R^{n+1}$ with holomorphic sections of $L^p$ with $L=\Oc(1)$ the tautological line bundle of $\P^{n+1}$. As in Example \ref{e:space}, we consider the standard Hermitian metrics on these line bundles. 
Note that $|z_1|^2+\cdots+|z_{n+1}|^2$ is constant on $\S^n$ and therefore, the formula for $\nu_p^\beta$ is simpler than the one in Example \ref{e:space}.
Observe also that a section of $L^p$ vanishes on $X$ if and only if it vanishes on $\S^n$. Therefore, $\{P_1,\ldots, P_{N_p}\}$ is identified to a basis of $H^0(X,L^p)$. We can now apply Theorem \ref{T:abstract} and Corollary \ref{C:manifold}.
\end{example}

The  plan  of the paper is  as  follows.  In Section \ref{S:regularity}, we discuss different notions of regularity for  the weighted compact set  $(K,\phi)$ and the measure $\mu$. We also give criteria to check the regularity conditions used in our study.  
In Section  \ref{S:almost-Fekete}, we prove the main theorem (Theorem  \ref{T:abstract}) which uses 
 an equidistribution result for almost Fekete configurations. 
The last  result has been obtained  in collaboration with Ma  in the last version of \cite[Remark 3.17]{DinhMaNguyen}.
For the reader's convenience, we provides here a  detailed proof  that we need in this paper.
Note that the case of Fekete points can be seen as the limit case of $\beta$-ensembles when $\beta\to\infty$. We refer to \cite{BBW, DinhMaNguyen, Leja57a, Leja57b, LevOrtega, Vu}, the references therein and also the end of this paper for more results on Fekete points and other configurations.

 \medskip\noindent
{\bf Acknowledgment.} 
The first author was  supported by Start-Up 
Grant R-146-000-204-133 from National University of Singapore. The  paper  was partially written during the  visits of the second author
at National  University of Singapore and Korea Institute for Advanced Study.  He would like 
to thank  these  organizations  for their financial support and their very warm hospitality.

\section{Pluri-regularity for weighted compact sets and measures}  \label{S:regularity}

As we have seen in Introduction, our study requires some regularity properties of the weighted compact set $(K,\phi)$ and  the probability measure $\mu$ on $K$. In this section, we will recall some known facts and also introduce and study  new notions that will be used in the proof of our main theorem. The reader will find in \cite{Demailly12, DinhSibony4, Hoermander, MaMarinescu07} basic notions and results from complex geometry and pluripotential theory.

Let $L$ be a positive (i.e., ample) holomorphic line bundle over 
a projective manifold $X$ of dimension $n$. Fix a smooth Hermitian metric $h_0$ on $L$ such that its first Chern 
form $\omega_0$ is a K\"ahler form on $X$. Define
$\mu^0:=\|\omega_0^n\|^{-1}\omega_0^n$ the probability measure 
associated with the volume form $\omega_0^n$. Here, $\|\omega_0^n\|$ is the total mass of $\omega_0^n$ which is the integral of this volume form on $X$.  
Recall that a real-valued function on $X$ is {\it quasi-p.s.h.} if it is locally the difference between a p.s.h. function and a smooth one. 
Let $\PSH(X,\omega_0)$ be the cone of {\it $\omega_0$-p.s.h. functions}, i.e., the quasi-p.s.h. functions 
$\varphi$ such that $\ddc\varphi+\omega_0\geq 0$. 

\begin{definition} \rm \label{def_weighted_K}
We call {\it weighted compact subset of $X$} a data $(K,\phi)$, 
where $K$ is a non-pluripolar compact subset of $X$ and 
$\phi$ is a real-valued continuous function on $K$. 
The function $\phi$ is called {\it a weight} on $K$. 
The {\it equilibrium weight} associated with $(K,\phi)$ is
 the upper semi-continuous regularization $\phi_K^*$ of 
 the function 
$$ \phi_K(z):=\sup\big\{\psi(z) : \psi \ \  \omega_0 \text{-p.s.h. such that } 
\ \psi\leq \phi \text{ on } K\big\}.$$
We also call {\it equilibrium measure of $(K,\phi)$} 
the normalized Monge-Amp\`ere measure 
$$\mu_\eq(K,\phi):=\|\omega_0^n\|^{-1}(\ddc\phi_K^*+\omega_0)^n.$$
\end{definition}

Note that the equilibrium measure $\mu_\eq(K,\phi)$ 
is a probability measure supported by $K$ and $\phi_K^*=\phi_K$ almost everywhere 
 with respect to this measure, see e.g., \cite{BermanBoucksom}. 
 The following notions are important in our study, see  \cite{DinhMaNguyen}.

\begin{definition} \rm \label{def_K_regular}
Denote by $P_K$ the projection onto $\PSH(X,\omega_0)$ 
which associates $\phi$ with $\phi_K^*$. We say that $(K,\phi)$
 is {\it regular} if $\phi_K$ is upper semi-continuous,
 i.e., $P_K\phi=\phi_K$. Let $(E, \|\ \|_E)$ be  a normed vector space of continuous functions on $K$ and $(F,\| \ \|_F)$ a normed vector space of functions on $X$. 
We say that $K$ is {\it $(E,F)$-regular}  
if $(K,\phi)$ is regular for  $\phi\in E$ 
and if the  projection $P_K$  sends bounded subsets of $E$ into bounded subsets of $F$. 
\end{definition}

We have the following result.

\begin{theorem}[\cite{DinhMaNguyen, Vu}] \label{t:DMNV}
Let $X$ and $L$ be as above. Let $K$ be a smooth compact real manifold in $X$ with piecewise smooth boundary. Assume that the tangent space to $K$ at any point is not contained in a complex hyperplane of the tangent space to $X$ at that point. Let $0<\alpha<1$ be any real number. Then $K$ is $(\Cc^\alpha,\Cc^{\alpha/2})$-regular. Moreover, it is $(\Cc^\alpha,\Cc^{\alpha})$-regular if the boundary of $K$ is smooth.
\end{theorem}

Consider now a real-valued function $\psi$ on $X$. 
We can associate the line bundle $L$ with a singular Hermitian metric
$h:=e^{-2\psi}h_0$. More precisely, if $v$ is a vector in the fiber of $L$ over 
a point $x\in X$, its norms with respect to the metrics $h$ 
and $h_0$ are  related by the formula
$$|v|_h=e^{-\psi(x)} |v|_{h_0}.$$
The metrics $h_0$ and $h$ induce in a canonical way metrics
 $h_0^{\otimes p}$ and $h^{\otimes p}$ on the power $L^{p}$ 
 of $L$.  They are related by the formula 
 $h^{\otimes p}=e^{-2p\psi} h_0^{\otimes p}$. 
Recall that for simplicity, we will use the notation $|\cdot|_{p\psi}$
 instead of $|\cdot|_{h^{\otimes p}}$ for the norm of a vector in
  $L^{p}$ with respect to the metric  $h^{\otimes p}$. We also drop the subscript $h_0$, e.g., $|v|$ means $|v|_{h_0}$.

Consider now a weighted compact set $(K,\phi)$ in $X$. We can, in a similar way, define the metric $h=e^{-2\phi}h_0$ on $L$ over $K$.  
Let  $\mu$ be a probability measure with support in $K$.  
Consider the natural $L^\infty$ and $L^2$ semi-norms on $H^0(X,L^p)$ 
 induced by the metric $h$ on $L$ and the measure $\mu$, which are defined for $s\in H^0(X,L^p)$ by
\begin{align}\label{eq:n1.15}
\|s\|_{L^\infty(K,p\phi)}:=\sup_K|s|_{p\phi} \qquad \text{and} \qquad 
\|s\|_{L^2(\mu,p\phi)}^2
:=\int_X |s|_{p\phi}^2 d\mu.
\end{align} 
We will only use measures $\mu$ such that the above semi-norms are  norms, i.e., there is no section $s\in H^0(X,L^p)\setminus \{0\}$ which vanishes on $K$ or on  the support of $\mu$. 
The first semi-norm is a norm when $K$ is not contained in a hypersurface of $X$. The second one is a norm
when the support of $\mu$   is not contained in a hypersurface of $X$. In particular, this is the case when $\mu$
is the normalized Monge-Amp\`ere measure with continuous potentials because such a measure has no mass on hypersurfaces of $X$. 

\smallskip

We need the following quantified Bernstein-Markov property, see also \cite{BermanBoucksom, BBW, Levenberg, NguyenZeriahi}.

\begin{definition}\label{D:BM} \rm
Let $\delta$ be a real number with $0<\delta<1$ and $(K,\phi)$ a weighted compact subset of $X$.
We say  that  a  positive  measure $\mu$ on $K$
   is  {\it  $\delta$-Bernstein-Markov} with respect to $(K,\phi)$   if
  there is a constant $A>0$  such that
  \begin{equation}\label{e:BM}
\|s\|_{L^\infty(K,p\phi)} \leq A e^{Ap^{1-\delta}} \|s\|_{L^2(\mu,p\phi)}\quad \text{for}\quad s\in H^0(X,L^p) \quad \text{and} \quad p\geq 1.
\end{equation}
If $\mu$ is  $\delta$-Bernstein-Markov  with respect to  $(K ,\phi)$ for all $0<\delta<1,$ then we say that
 $\mu$ is  {\it $1$-Bernstein-Markov} with respect to $(K,\phi)$.
\end{definition}

The following lemma shows that we can use the notion for other norms $L^r$. 

\begin{lemma} \label{l:Lr}
Let $\delta, r$ be real numbers with $0<\delta<1$ and $r >0$. Let  $(K,\phi)$ be a weighted compact subset of $X$ and $\mu$  a positive measure on $K$. Then $\mu$ is $\delta$-Bernstein-Markov with respect to $(K,\phi)$ if and only if 
  there is a constant $A'>0$  such that
$$\|s\|_{L^\infty(K,p\phi)} \leq A' e^{A'p^{1-\delta}} \|s\|_{L^r(\mu,p\phi)}\quad \text{for}\quad s\in H^0(X,L^p) \quad \text{and} \quad p\geq 1.$$
\end{lemma}
\proof
Assume that $\mu$ is $\delta$-Bernstein-Markov with respect to $(K,\phi)$. We will only show the existence of $A'$ as in the lemma because the converse property can be  obtained in the same way. 
So we have property \eqref{e:BM}. 
Without loss of generality, we can assume that $\mu$ is a probability measure.
If $r\geq 2$, then the $L^r$-norm is larger or equal to the $L^2$-norm. Therefore, we can just take $A':=A$.

Assume now that $0<r<2$. 
By H\"older's inequality, we have 
$$\|s\|_{L^2(\mu,p\phi)} \leq \|s\|_{L^r(\mu,p\phi)}^{r/2} \|s\|_{L^\infty(K,p\phi)}^{1-r/2}.$$
This, together with  \eqref{e:BM}, gives us the desired property for a suitable value of $A'$. 
\endproof

In order to get a simple criterium for a measure to have the  $\delta$-Bernstein-Markov property, we need the following notion.

\begin{definition}\label{D:max_set} \rm
A compact set $W$  is  said to satisfy {\it the maximum principle} relatively to  a  weighted compact set $(K,\phi)$ 
if $W\subset K$ and  
$$
\sup_{K} (\psi-\phi)=\sup_W (\psi-\phi) \quad\text{for every}\quad \psi\in\PSH(X,\omega_0).
$$
\end{definition}
Clearly, $W=K$   satisfies  the maximum principle relatively to    $(K,\phi).$
In general, $W$  may be much  smaller  than  $K$, see Remark \ref{rk_max_Pn} below.

\begin{proposition}\label{P:max_set}
Let $(K,\phi)$ be a   weighted compact set and $W$ a compact subset of $K$. Define
$$\partial^\phi_{\omega_0} K:=\left\lbrace   z\in K: P_K\phi(z)=\phi(z)\right\rbrace.$$
Then $W$ satisfies the maximum principle relatively to   $(K,\phi)$ if and only if $W\cap \partial^\phi_{\omega_0} K$ satisfies the same property. 
In particular, $\partial^\phi_{\omega_0} K$
satisfies the maximum principle relatively to   $(K,\phi)$
\end{proposition}
\proof
Observe that the second assertion is a consequence of the first one and Definition \ref{D:max_set} by taking $W=K$. We prove now the first assertion. 
If  $W\cap \partial^\phi_{\omega_0} K$ satisfies the maximum principle relatively to $(K,\phi)$, then clearly $W$ satisfies the same property. Assume that 
 $W$ satisfies this maximum principle. It remains to prove the same property for $W\cap \partial^\phi_{\omega_0} K$.
 
 Recall that $P_K\phi$ is upper semi-continuous and $\phi$ is continuous. Since $P_K\phi\leq \phi$, we deduce that 
 $$\partial^\phi_{\omega_0} K=\left\lbrace   z\in K: P_K\phi(z)\geq \phi(z)\right\rbrace$$
 So it is a compact set. 
 
Let $ \psi\in\PSH(X,\omega_0)$ and set $m:= \max_{K} (\psi-\phi).$ Note that $\psi$ is also upper semi-continuous. 
Since $W$ satisfies the maximum principle relatively to $(K,\phi)$, there is a point $z_0\in W$ such that $\psi-\phi$ attains its maximum value  at $z_0.$
We have $\psi(z_0)-m=\phi(z_0)$ and 
 $\psi-m\leq \phi$ on $K.$ Since  $ \psi-m\in\PSH(X,\omega_0),$ the last inequality implies that $\psi-m\leq P_K\phi$.
In particular,
$$\psi(z_0)-m\leq (P_K\phi)(z_0).$$
This, combined  with  the equality  $\psi(z_0)-m=\phi(z_0)$ and  the inequality $(P_K\phi)(z_0)\leq \phi(z_0),$ implies that $(P_K\phi)(z_0)= \phi(z_0).$ Hence, $z_0\in \partial^\phi_{\omega_0} K$ and the proposition follows.
\endproof

\begin{remark}\rm 
By \cite[Prop. 2.10, Cor. 2.5]{BermanBoucksom}, the equilibrium measure $\mu_\eq(K,\phi)$ is  supported by $\partial^\phi_{\omega_0} K$ and its support also satisfies the maximum principle. 
\end{remark}

\begin{remark}\rm \label{rk_max_Pn}
Let $X$ be the projective space $\P^n$, seen as 
 the natural compactification of $\C^n$. Let  $L$ be the tautological
 line bundle $\mathcal O(1)$ over $\P^n$. Then the holomorphic sections of  $L^{p}=\mathcal O(p)$ can be identified to the
 complex polynomials of degree $\leq p$ on $\C^n$. With the standard Fubini-Study metric on $O(p)$, if a section $s$ of $L^p$ corresponds to a polynomial $P(z)$ of degree $\leq p$, then 
 $$|s(z)|=|P(z)|(1+\|z\|^2)^{-p/2}.$$
 Consider a compact subset $K$ of $\C^n$ and take  $\phi:=-{1\over 2} \log (1+\|z\|^2)$ on $K.$ It is not difficult to check that the boundary of $K$ satisfies  the maximum principle relatively to $(K,\phi).$ 
\end{remark}

\begin{theorem}\label{T:Bernstein-Markov}
Let $X,L,h_0$   be as above, $(K,\phi)$ a weighted compact 
subset of $X$ and $\mu$  a probability measure on $K$. Let $W\subset K$ be a  compact set and $0<\delta<1$  a real number. 
Assume in addition  the  following  conditions:
\begin{enumerate}
\item[(i)] the functions $\phi$ and $P_k\phi$ are H\"older continuous;
\item[(ii)]  $W$ satisfies  the maximum principle relatively to $(K,\phi);$
\item[(iii)]  $\mu$ satisfies  the mass-density condition  with respect to $W$, see Definition \ref{D:mass-density}.
\end{enumerate}
Then  $\mu$ is  a $\delta$-Bernstein-Markov measure with respect to $(K,\phi).$ 
\end{theorem}

\begin{remark}\rm \label{rk_BM}
We will see in the proof of this theorem that the condition (i) can be replaced by the following much weaker condition :
there are constant $c>0$  such that for  $z\in X$ and $w\in K$
$$ |(P_K \phi)(z)-(P_K\phi)(w)| \leq c (1+\log^-\dist(z,w))^{-{\delta\over 1-\delta}},$$
and for $z,w\in K$
$$ |\phi(z)-\phi(w)|\leq  c(1+\log^-\dist(z,w))^{-{1\over 1-\delta}},$$
where $\log^-:=\max(0,-\log)$.
\end{remark}

We are inspired by an idea of Bloom \cite[Theorem 4.1]{Bloom}.
Define  $\epsilon:=p^{-\delta}$  and $r:=e^{-c'p^{1-\delta}}$ where $c'>0$ is  a large enough constant independent of $p$.
It follows from Assumption (i) (see also Remark \ref{rk_BM}) that
\begin{equation}\label{e:P_K-Hoelder}
\big|(P_K  \phi\big)(z)-  (P_K\phi)(z_0) \big|\leq \epsilon \quad\text{for}\quad z\in B(z_0,2r) \quad \text{and} \quad z_0\in K.
\end{equation}

Fix $p\geq 1$ and $s\in H^0(X, L^p)\setminus\{0\}$. We need to prove  inequality \eqref{e:BM} for some constant $A>0$ independent of $p$ and $s$.
Observe  that
$$\ddc  {1\over p} \log|s|={1\over p} [s=0] -\omega_0\geq -\omega_0,$$
where $[s=0]$ is the current of integration on the hypersurface $\{s=0\}$. 
So 
 ${1\over p} \log|s|$ is $\omega_0$-p.s.h. 
This,  together  with Assumption (ii), implies
the  existence  of a  point
 $z_0\in W$  such that 
\begin{equation}\label{e:z_0}
|s(z_0)|_{p\phi}=\max_{z\in K} |s(z)|_{p\phi}.
\end{equation}

\begin{lemma}\label{l:s_for_z_near_z0_bis}
We have 
$$\big| |s(z)|-|s(z_0)|\big |\leq {1\over 4}    |s(z_0)| \quad\text{for}\quad z\in  B(z_0, r^2) .$$
\end{lemma}
\proof
Consider a section
$s'=cs$ where the constant $c$ is chosen so that $ \|s'\|_{L^\infty(K,p\phi )}=1$. The last property implies the inequality
 ${1\over p} \log|s'|\leq \phi$ on $K.$
We  have  seen that
  ${1\over p} \log|s|$ is $\omega_0$-p.s.h. So $s'$ satisfies a similar property. Hence,
 ${1\over p} \log|s'|\leq P_K\phi$ on $X$. We then deduce  the following
 Bernstein-Walsh type inequality  
 $$|s(z)|\leq  \|s\|_{L^\infty(K,p\phi )}e^{p (P_K\phi)(z)}\quad \text{for} \quad z\in X.$$
Using \eqref{e:z_0}, we rewrite the last  inequality   for $z\in B(z_0,2r)$ as 
 \begin{equation*}
 |s(z)|\leq  |s(z_0)|_{p\phi} e^{p (P_K\phi)(z)}= |s(z_0)| e^{-p\big(\phi(z_0)- (P_K\phi)(z_0)\big)}
 e^{p\big((P_K\phi)(z)- (P_K\phi)(z_0)\big)}.
 \end{equation*}
 Using $\phi(z_0)\geq (P_K\phi)(z_0)$ and   \eqref{e:P_K-Hoelder}, we obtain
\begin{equation}\label{e:s_for_z_near_z0}
|s(z)|\leq  |s(z_0)| e^{p \epsilon}\quad \text{for} \quad z\in B(z_0,2r).
\end{equation}

Let $\sigma$ be a holomorphic frame for $L$ on an open neighborhood
$U$ of   $z_0 $  with $|\sigma(z_0)|=1.$ Write  $s=h\sigma^{\otimes p}$ with $h$ a holomorphic function on $U$. Using local  coordinates near $z_0$ and shrinking
$U$ if necessary,  we  may identify $U$ with the open unit ball in $\C^n.$  
We can also assume that 
$$\big| |\sigma(z)| - 1 \big| \leq c \|z-z_0\|$$
for some constant $c>0$ independent of $z\in U$. 
For $z\in B(z_0,2r)$, we have $\|z-z_0\|\ll p^{-1}$ and the previous inequality 
 implies that $|\sigma^{\otimes p}(z)|$ belongs to the interval $[7/8,9/8]$  when $z\in B(z_0,2r)$. So the norm $|s(z)|$ is bounded below and above by $7|h(z)|/8$ and $9|h(z)|/8$ respectively. 

Consider the unit  vector  $
v:={z-z_0\over \| z-z_0\|}$ in $\C^n,$ and   the following holomorphic function of one  variable
$$f(\zeta):=h(z_0+\zeta v),\qquad \zeta\in\D.$$
We have for $z\in B(z_0,r)$
\begin{equation}\label{e:h}
|h(z) -h(z_0)|=|f(\|z-z_0\|)-f(0)|\leq  \|z-z_0\| \sup_{|\zeta|\leq r}|f'(\zeta)|.
\end{equation}
On the other hand, for $|\zeta|\leq 2r,$ we have $(z_0+\zeta v)\in B(z_0,2r),$ and 
by using the definition of $f,$ $h$,   \eqref{e:s_for_z_near_z0} and that $|s(z)|$ is in-between $7|h(z)|/8$ and $9|h(z)|/8$,  we obtain
\begin{equation*}
\sup_{|\zeta|\leq 2r}|f(\zeta)|\leq  c|s(z_0)| e^{p \epsilon}  
\end{equation*}
for some constant $c>0$. By Cauchy's formula,
$$\sup_{|\zeta|\leq r}|f'(\zeta)|\leq {c\over r}  |s(z_0)| e^{p \epsilon}.$$
This, together with \eqref{e:h} and the choice of $\epsilon, r$, implies for $z\in B(z_0,r^2)$ that 
\begin{equation*} 
 |h(z) -h(z_0)|\leq  cr |s(z_0)| e^{p \epsilon}
 \ll  |s(z_0)| .
 \end{equation*}
Recall that $|h(z_0)|=|s(z_0)|$ and $|s(z)|$ is bounded by $7|h(z)|/8$ and $9|h(z)|/8$. So the last inequality implies the lemma.
\endproof

\noindent
{\bf End of the proof of Theorem \ref{T:Bernstein-Markov}.} 
We only need to consider $p$ large enough.
We will prove that
\begin{equation}\label{e:s_for_z_near_z0_weight}
|s(z)|_{p\phi}\geq {1\over 2}\|s\|_{L^\infty(K,p\phi )}\quad\text{for}\quad z\in 
K\cap B(z_0, r^2) .
\end{equation}
We have
$$
\big| |s(z)|_{p\phi}-|s(z_0)|_{p\phi} \big|\leq \big| |s(z)|_{p\phi(z_0)}-|s(z_0)|_{p\phi(z_0)}\big|
+ \big| |s(z)|_{p\phi(z)}-|s(z)|_{p\phi(z_0)}\big|.
$$
Denote respectively by $A_1$ and $A_2$ the first and second terms in the last sum. 
By Lemma \ref{l:s_for_z_near_z0_bis}, we  have
$$A_1 \leq  {1\over 4}  |s(z_0)| e^{-p\phi(z_0)}.$$
On the  other  hand,  by Lemma \ref{l:s_for_z_near_z0_bis} again, we have
$$ A_2=|s(z)||e^{-p\phi(z)}-e^{-p\phi(z_0)}|\leq 2 |s(z_0)|e^{-p\phi(z_0)}|1-e^{-p(\phi(z)-\phi(z_0))}|.$$
 Since  $z\in B(z_0, r^2),$  we deduce  from   Assumption (i) of the theorem  (see also Remark \ref{rk_BM})  that $|p(\phi(z)-\phi(z_0))|\leq 1/16.$
 Hence,
 $|1-e^{-p(\phi(z)-\phi(z_0))}|\leq 1/8.$
 Combining  the above estimates for $A_1$ and $A_2,$ we obtain
 $$
 \big||s(z)|_{p\phi}-|s(z_0)|_{p\phi}\big|\leq {1\over 2}|s(z_0)|_{p\phi}\quad\text{for}\quad z\in  
K\cap B(z_0, r^2) .
 $$
 This, combined with \eqref{e:z_0}, implies \eqref{e:s_for_z_near_z0_weight}.
 
 Now, using \eqref{e:s_for_z_near_z0_weight} and  Assumption (iii), we get 
\begin{eqnarray*}
\int_K |s(z)|^2_{p\phi} d\mu &\geq &  \int_{K\cap B(z_0,r^2)} |s(z)|^2_{p\phi} d\mu\\
&\geq &\Big( \min\limits_{K\cap B(z_0,r^2)}  |s(z)|^2_{p\phi}\Big)
\mu( K\cap B(z_0,r^2))\\
&\geq& {1\over 4} c  r^{2\rho}\|s\|^2_{L^\infty(K,p\phi)},
\end{eqnarray*}
where $c>0$ is the constant in Definition \ref{D:mass-density}.
Hence,
$$
\|s\|_{L^\infty(K,p\phi)} \leq 2c^{-1/2} e^{\rho c' p^{1-\delta}}\|s\|_{L^2(\mu,p\phi)}.
$$
So $\mu$ is $\delta$-Bernstein-Markov  with respect to  $(K, \phi).$ 
\hfill $\square$

\medskip

We have the following result where Condition (ii) is automatically satisfied for $W=K$. It allows us to obtain Corollary \ref{C:manifold} as a direct consequence of Theorem \ref{T:abstract}. Note that in Corollary \ref{C:manifold} we only need to assume that the measure $\mu$ satisfies the mass-density condition with respect to a compact $W\subset K$ which satisfies the maximum principle relatively to $(K,\phi)$. 

\begin{corollary}\label{C:Bernstein-Markov}
Let $X,L,h_0$ be as above, $K$  a  compact 
subset of $X,$ $W$ a compact subset of $K$ and $\mu$ a probability  measure on $K.$ 
Assume in addition  the  following  conditions:
\begin{enumerate}
\item[(i)] $K$ is $(\Cc^\alpha,\Cc^{\alpha'})$-regular for some constants $\alpha>0$ and $\alpha'>0$; 
\item[(ii)] $W$  satisfies the maximum principle  relatively to $(K,\phi);$
\item[(iii)]  $\mu$ satisfies  the mass-density condition  with respect to $W.$
\end{enumerate}
Then $\mu$ is  a $1$-Bernstein-Markov measure with respect  to  $(K,\phi)$   for every  $\phi\in\Cc^\alpha(K).$  
\end{corollary}
\proof
Since  $\phi\in\Cc^\alpha(K)$ and $K$ is $(\Cc^\alpha,\Cc^{\alpha'})$-regular,
$(K,\phi) $ satisfies the hypotheses of Theorem \ref{T:Bernstein-Markov}. According to that theorem,  $\mu$ 
is $\delta$-Bernstein-Markov with respect to $(K,\phi)$  for every $0<\delta<1.$ The corollary follows.
\endproof

\section{Almost-Fekete configurations and proof of the main result}  \label{S:almost-Fekete}

In this section, we will give the proof of the main theorem. An important ingredient is the 
equidistribution of almost-Fekete points towards the equilibrium measure. This property is already mentioned in the last version of \cite{DinhMaNguyen}, see also \cite{Levenberg2}. For the reader's convenience, we will give here some details. We also give at the end of this section another application of this result.

\begin{theorem}[\cite{DinhMaNguyen}] \label{thm_DMN}
Let $X,L,h_0$ be as above and $K$ a compact 
subset of $X$. Let $0<\alpha\leq 2,$ $ 0<\alpha'\leq 1$ and $0<\gamma\leq 2$ be constants. 
Assume that $K$ is $(\Cc^\alpha,\Cc^{\alpha'})$-regular. Let $\phi$ be a $\Cc^\alpha$ real-valued function on $K$ 
 and $\mu_\eq(K,\phi)$ the equilibrium measure associated 
 with the weighted set $(K,\phi)$. Then, there is a constant $c>0$  with the following property.
 For every $p\geq 1$ and
every configuration
$\x=(x_1,\ldots,x_{N_p})\in K^{N_p},$ denote by $\mu^{\x}$ the empirical measure associated with $\x$ and   let 
$S_p$ be any basis of $H^0(X,L^p)$. 
Define 
$$\sigma_{\x}:={1\over pN_p} \log \|\det S_p\|_{L^\infty(K,p\phi)}-{1\over pN_p} \log \|\det S_p(\x)\|_{p\phi}.$$
Then   we have for all $p>1$
$$\dist_\gamma(\mu^{\x}, \mu_\eq(K,\phi))\leq c p^{-\alpha''\gamma}(\log p)^{3\alpha''\gamma}+c \sigma_{\x}^{\gamma/4}\quad \text{with} 
\quad\alpha'':=  \alpha'/(24+12\alpha').$$
 \end{theorem}

Note that $\det S_p$ is a section of the line bundle $(L^p)^{\boxtimes N_p}$ over $X^{N_p}$. The given metric $h_0$ on $L$ and 
the weight $\phi$ induces naturally a metric and a weight for this line bundle. So  $\|\det S_p\|_{L^\infty(K,p\phi)}$ is the sup-norm of $\det S_p$ on $K^{N_p}$ and $\|\det S_p(\x)\|_{p\phi}$ is the norm of the value of this section at the point $\x$. Both of them are measured
using the above natural metric and weight. Observe that $\sigma_\x$ is independent of the choice of $S_p$ and we always have $\sigma_\x\geq 0$. When $\sigma_\x=0$, the point $\x$ is called {\it a Fekete configuration of order $p$} of $L$ with respect to the weighted compact set $(K,\phi)$. The theorem shows that if $\sigma_\x$ is small enough (e.g., when $\sigma_\x=0$),  then $\mu^\x$ tends to $\mu_\eq(K,\phi)$ as $p\to\infty$. 

\smallskip

We now sketch the proof of Theorem \ref{thm_DMN}.
Recall that the {\it Monge-Amp\`ere energy functional} 
$\mE $, defined on bounded weights in 
$\PSH(X,\omega_0)$, is characterized by
$$\left. {d\over dt}\right|_{t=0} \mE ((1-t)\phi_1+t\phi_2)
=\|\omega_0^n\|^{-1} \int_X(\phi_2-\phi_1)  (\ddc\phi_1+\omega_0)^n.$$
So $\mE $ is only defined up to an additive constant, but the differences such as 
$\mE(\phi_1) -  \mE(\phi_2)$
are well-defined,  see \cite{BermanBoucksom}  and also \eqref{eq_normalization} below. 

Consider a non-pluripolar compact set $K\subset X$ and a continuous weight $\phi$ on $K$. Define the {\it energy at the equilibrium weight} of $(K,\phi)$ as
$$\mE_\eq(K,\phi):= \mE(P_K\phi).$$
This functional  is also  well-defined up to an additive constant. We have the following property.

\begin{lemma}[\cite{BermanBoucksom}, Th. B]\label{lem_differentiability}
The map  $\phi\mapsto \mE_\eq(K,\phi),$  
defined on the affine space of continuous weights on $K$,
is concave  and G\^ateaux differentiable, with directional derivatives 
given by integration against the equilibrium measure:
$$\left . {d\over dt}\right|_{t=0}  \mE_\eq (K,\phi+tv)
=\big\langle v ,\mu_\eq(K,\phi)  \big\rangle \quad
\text{for every continuous function } v \text{ on } K.$$ 
\end{lemma}

Let $\mu$ be a probability measure on $X$ and $\phi$ a continuous function on the support of $\mu$. The semi-norm
 $\|\cdot\|_{L^2(\mu,p\phi)}$ on $H^0(X,L^{p})$ is defined as in \eqref{eq:n1.15} and recall that we only consider measures $\mu$ for which this semi-norm is a norm.  
Let $\mathcal B^2_p(\mu,\phi)$ denote the unit ball in  
$H^0(X,L^{p})$ with respect to this norm and $N_p:=\dim H^0(X,L^{p})$.
Consider the following 
$\mL_p$-functional
\begin{align}\label{eq:n4.5}
\mL_p(\mu,\phi):=  {1\over  2pN_p}  
\log\vol\mathcal B^2_p(\mu,\phi). 
\end{align}
Here, $\vol$ denotes the Lebesgue measure on the 
vector space $H^0(X,L^{p})$ which depends on the choice of an Euclidean norm on   $H^0(X,L^{p})$. So the volume is  only defined up 
to a multiplicative constant. Nevertheless, 
the  differences such as  
$\mL_p(\mu_1,\phi_1)- \mL_p(\mu_2,\phi_2)$ are well-defined 
and do  not depend on the choice of $\vol$ for any
 probability measures $\mu_1$ and $\mu_2$,  see \cite {BermanBoucksom} and also \eqref{eq_normalization} below.

Consider the norm
 $\|\cdot\|_{L^\infty(K,p\phi)}$ on $H^0(X,L^{p})$ defined   in \eqref{eq:n1.15}. 
Let $\mathcal B^\infty_p(K,\phi)$ denote the unit ball in 
$H^0(X,L^{p})$ with respect to this norm. Define 
$$  \mL_p(K,\phi):=  {1\over  2pN_p}  \log\vol\mathcal B^\infty_p(K,\phi).$$

Let  $\{s_1,\ldots, s_{N_p}\}$ be an orthonormal basis of
$H^0(X,L^{p})$ with respect to the above $L^2$-norm, see \eqref{eq:n1.15}.

\begin{definition} \rm \label{defi_Bergman_distortion} 
We call   {\it Bergman function of $L^{p}$}, 
associated with $(\mu,\phi)$, the function $\rho_p(\mu,\phi)$  on the support of $\mu$  given by
$$\rho_p(\mu,\phi)(x):=\sup\Big\{|s(x)|_{p\phi}^2: \ \  s\in H^0(X,L^p), \|s\|_{L^2(\mu,p\phi)}=1\Big\} =\sum_{j=1}^{N_p}|s_{j}(x)|_{p\phi}^2$$
and we define {\it the Bergman measure} associated with  $(\mu,\phi)$ by
$$\Bc_p(\mu,\phi):=N_p^{-1}\rho_p(\mu,\phi)\mu.$$
\end{definition}

It is not difficult to obtain the identity in the definition of $\rho_p(\mu,\phi)$ and  to check that $\Bc_p(\mu,\phi)$ is a probability measure. 
 Note also that when $\mu$ is the average of $N_p$ Dirac masses at generic points, one can easily deduce 
from Definition \ref{defi_Bergman_distortion} that $\Bc_p(\mu,\phi)=\mu$, by considering sections vanishing on $\supp(\mu)$ except at a point. Such sections exist because $N_p=\dim H^0(X,L^p)$. In fact, this property holds for all points $x_1,\ldots, x_{N_p}$ such that the section $\det S_p$ considered in Introduction does not vanish at $(x_1,\ldots, x_{N_p})$.

\begin{lemma} \label{lem_mL_p}
\begin{enumerate}
\item[(a)] The  functional  $\phi\mapsto \mL_p(\mu, \phi)$ is concave on
 the space of all   continuous  weights on the support of $\mu$.
\item[(b)] The directional derivative of $\mL_p(\mu,\cdot)$ at 
a continuous weight $\phi$ on the support of $\mu$ is given by the integration against
the Bergman measure $\Bc_p(\mu,\phi),$ that is,
$$\left.{d\over  dt}  \mL_p(\mu, \phi+tv)\right |_{t=0}
= \langle v,\Bc_p(\mu,\phi)\rangle ,\quad  \text{with } v,\phi \text{ continuous on the support of } \mu.$$
\item[(c)]   
Let $\mu$ be a probability measure with $\supp(\mu)\subset K$ such that the $L^2$-semi-norm in \eqref{eq:n1.15} is a norm. Assume also that $(K,\phi)$ is a regular weighted compact set. Then 
$$\mL_p(K,\phi)=\mL_p(X,P_K\phi)\quad\text{and}\quad \mL_p(K,\phi)\leq   \mL_p(\mu,\phi).$$
\end{enumerate}
 \end{lemma}
\proof
The   concavity property of the functional 
 $\mL_p$ in Part (a) has been established in 
  \cite[Proposition 2.4]{BBW}.
 Part (b) has been established 
in \cite[Lemma 5.1]{BermanBoucksom}. The property was stated there for smooth $\phi$ but the proof also works for continuous functions, see also \cite[Lemma 3.1]{Berndtsson2} and \cite[Lemma 2]{Donaldson}.
For Part (c), see \cite[Proposition 2.5, Lemma 3.4]{DinhMaNguyen}. 
\endproof

  From now on, in order to simplify the notation, we use the following normalization
\begin{equation}\label{eq_normalization}
 \mE_\eq(X,0)=0,\quad \mL_p(X,0)=0 \quad\text{and}
 \quad 
\mL_p(\mu^0,0)=0\quad \text{for} \quad p\geq 1.
 \end{equation}
Here, the function identically 0 is used as a smooth strictly $\omega_0$-p.s.h. weight. Recall also that
$\mu^0=\|\omega_0^n\|^{-1}\omega_0^n$ is the probability measure 
associated with the volume form $\omega_0^n$. 

\medskip

The  following result is  an immediate consequence of \cite[Proposition 3.10]{DinhMaNguyen}. Recall that $\Cc^{k,\alpha}=\Cc^{k+\alpha}$ for $0\leq \alpha<1$ and $\Cc^{k,1}$ is the space of $\Cc^k$ functions whose partial derivatives of order $k$ are Lipschitz.

\begin{proposition}\label{prop_L_to_E_4} 
Let $0<\alpha\leq 1$ and $A>0$ be constants.
Let $\phi $ be  an $\omega_0$-p.s.h. weight of class 
$\Cc^{0,\alpha}$  on $X$ such that $ \|\phi\|_{\Cc^{0,\alpha}} \leq A$.   Then, there is a constant
$c_{A,\alpha}>0$ depending only on $X,L,\omega_0, A$ and $\alpha$
such that we have for all  $p>1$
$$ \big | \mL_p(\mu^0,\phi) - \mE_\eq(X,\phi)\big| \leq  c_{A,\alpha}  (\log p)^{3\beta_\alpha } p^{-\beta_\alpha}$$
 and
$$\big | \big(\mL_p(X,\phi) -  \mE_\eq(X,\phi)\big| \leq  c_{A,\alpha}   (\log p)^{3\beta_\alpha} p^{-\beta_\alpha},$$
where $\beta_\alpha:= \alpha/(6+3\alpha)$.
\end{proposition}

For the following proposition, we refer to the discussion after Theorem \ref{thm_DMN} for the notation.

\begin{proposition} \label{prop_D_L_E} 
Let $K$ be a compact 
subset of $X$.
Let   $0<\alpha\leq 2$ and $ 0<\alpha'\leq 1$   be constants. 
Assume that $(K,\phi)$ is a weighted compact set  with $\phi\in  \Cc^\alpha(K)$
such that $K$ is $(\Cc^\alpha,\Cc^{\alpha'})$-regular.  Then there is a constant $c>0$ with the following property. 
For  $p\geq 1$ and $\x=(x_1,\ldots,x_{N_p})\in K^{N_p},$ denote by $\mu^{\x}$ the empirical measure associated with $\x$ and 
let $S_p$ be  a  basis of $H^0(X,L^p)$.    Define 
$$\sigma_{\x}:={1\over pN_p} \log \|\det S_p\|_{L^\infty(K,p\phi)}-{1\over pN_p} \log \|\det S_p(\x)\|_{p\phi}.$$
We have for all $p>1$
$$|\mL_p(\mu^{\x},\phi)-\mE_\eq(K,\phi)|\leq c \big (p^{-1}\log p+\sigma_{\x}+\big |\mL_p(\mu^0,P_K\phi)
-\mE_\eq(K,\phi)\big|\big).$$
\end{proposition}
 \proof
 Observe that $\sigma_\x$ does not depend on the choice of $S_p$. So choose $S_p$ which is an orthonormal basis of $H^0(X,L^p)$
 with respect to the $L^2$ -norm without weight. 
 Let $\mu_p$ be   the  empirical measure associated with a Fekete configuration of  order $p.$
  Using  identity  \cite[(2.4)]  {BBW}, we  get 
$$ {1\over  2p N_p}\log{  \vol \mathcal B^2_p(\mu^0,0) 
\over \vol\mathcal B^2_p(\mu_p,\phi)} =  {1\over  pN_p}
\log \| \det S_p\|_{L^\infty(K,p\phi)}        -{1\over 2p} \log N_p$$
and
$${1\over  2p N_p}\log{  \vol \mathcal B^2_p(\mu^0,0) 
\over \vol\mathcal B^2_p(\mu^{\x},\phi)} = {1\over pN_p} \log \|\det S_p(\x)\|_{p\phi} -{1\over 2p} \log N_p.$$
Subtracting the last line  from the previous   one and  using  \eqref{eq:n4.5}, we obtain
 $$\sigma_{\x}= \mL_p(\mu^{\x},\phi)-\mL_p(\mu_p,\phi).$$
On the  other hand,
with the normalization \eqref{eq_normalization},  \cite[Proposition 3.12]{DinhMaNguyen} tells us that there is a constant $c>0$ satisfying 
$$|\mL_p(\mu_p,\phi)-\mE_\eq(K,\phi)|\leq c \big (p^{-1}\log p +\big |\mL_p(\mu^0,P_K\phi)
-\mE_\eq(K,\phi)\big|\big)\quad\text{for}\quad p>1.$$
This, combined with the previous identity, implies  the proposition.
 \endproof
 
The following two lemmas were obtained in \cite[Lemmas 3.13 and 3.14]{DinhMaNguyen}.

\begin{lemma}\label{lem_MA_Lip} There is  a constant $c>0$ such that 
for every  continuous weight  $\phi$ on $K$  and every function $v$ of class $\Cc^{1,1}$ on $X$, we have   
$$\big|\langle \mu_\eq(K, \phi+tv) -\mu_\eq(K, \phi),v\rangle\big|\leq  c|t|\|v\|_{L^\infty(K)}\|\ddc v\|_\infty \quad \text{for}\quad t\in\R.$$
\end{lemma}

\begin{lemma}\label{key_lemma}
Let $ \epsilon>0$ and $M>0$ be constants. Let  $F$ and $G$ be   functions defined on $[-\epsilon^{1/2},\epsilon^{1/2}]$  such that
\begin{enumerate}
\item[(i)] $F(t)\geq G(t)-\epsilon$ and 
$|F(0) -  G(0)|\leq \epsilon;$
\item[(ii)]  $F$ is  concave on $[-\epsilon^{1/2},\epsilon^{1/2}]$ and differentiable at $0;$
\item[(iii)]  $G$ is differentiable in $ [-\epsilon^{1/2},\epsilon^{1/2}],$ and 
its derivative $G'$  satisfies
  $|G'(t)-G'(0)|\leq M\epsilon^{1/2}$ for $t \in [-\epsilon^{1/2},\epsilon^{1/2}] .$ The last inequality holds when $|G'(t)-G'(0)|\leq M|t|$. 
\end{enumerate}
Then we have
$$|F'(0)- G'(0)|\leq (2+M)\epsilon^{1/2}.$$
\end{lemma}

\medskip
 \noindent {\bf End of the proof of Theorem \ref{thm_DMN}.}
By \eqref{eq:n1.23}, we only need to consider the case $\gamma=2,$ i.e., to prove that
 \begin{equation*}
   \big |  \langle \mu^{\x} -\mu_\eq(K,\phi), v \rangle \big | 
   \lesssim    p^{-2\alpha''} (\log p)^{6\alpha''}  +\sigma_\x^{1/2}
 \end{equation*}
 for every test $\Cc^2$ function $v$ such that
 $\| v\|_{\Cc^2}\leq 1.$ Recall that $\alpha'':=\alpha'/(24+12\alpha')$.

Define 
$$F(t):=\mL_p(\mu^{\x},\phi+tv)\qquad \text{and}\qquad  G(t)
:= \mE_\eq(K,\phi+tv)=\mE_\eq(X,P_K(\phi+tv))$$
for $t$ in a neighborhood of $0\in\R$. 
By  Lemma \ref{lem_mL_p}(c),
 $$\mL_p(\mu^\x,\phi+tv)\geq  \mL_p(K,\phi+tv)= \mL_p(X,P_K(\phi+tv)).$$
As $0<\alpha\leq 2,$ we infer 
$\phi+tv\in\Cc^\alpha(K).$
Since $K$ is $(\Cc^\alpha,\Cc^{\alpha'})$-regular,  we deduce that $P_K(\phi+tv)$ is an 
$\omega_0$-p.s.h. weight on $X$ with bounded $\Cc^{\alpha'}$-norm.
This, coupled with the second inequality in 
Proposition \ref{prop_L_to_E_4}, applied to $P_K(\phi+tv)$  and $\alpha'$ instead of $\alpha$,  shows that
 \begin{align}\label{eq:n4.91} F(t)- G(t) \gtrsim -   p^{-4\alpha''} (\log p)^{12\alpha''}.
 \end{align}
 
An application of the first inequality in  Proposition \ref{prop_L_to_E_4} for $\alpha'$ instead of $\alpha$   gives 
$$     \big |\mL_p(\mu^0,P_K\phi) -\mE_\eq(K,\phi)\big|  \lesssim p^{-4\alpha''} (\log p)^{12\alpha''}.$$ 
Consequently, applying Proposition \ref{prop_D_L_E} yields   
$$|F(0)-G(0)| \lesssim    p^{-4\alpha''} (\log p)^{12\alpha''}+\sigma_{\x}.$$
Recall from Lemma \ref{lem_mL_p}(a) that $F$ is concave.
Moreover, by Lemma  \ref{lem_mL_p}(b), we have
$$F'(0)=\langle v, \Bc_p(\mu^{\x},\phi)  \rangle.$$
On the other hand, by Lemma \ref{lem_differentiability}, $G$
 is differentiable with
\begin{align}\label{eq:n4.94}
G'(t)=\langle v, \mu_\eq( K,\phi+tv)  \rangle.
\end{align}

Finally, by  Lemma \ref{lem_MA_Lip}, condition (iii) in  
Lemma \ref{key_lemma} is satisfied for a suitable constant $M>0$. 
Combining this and the discussion between \eqref{eq:n4.91}-\eqref{eq:n4.94}, we are in the position to apply Lemma \ref{key_lemma} to a constant $\epsilon$ of order  $p^{-4\alpha''} (\log p)^{12\alpha''}+\sigma_{\x}$.
Using  the  above expression for $F'(0)$ and $G'(0),$  we get 
$$\big | \langle \Bc_p(\mu^{\x},\phi) , v \rangle
 - \langle \mu_\eq(K, \phi), v  \rangle \big |  \lesssim p^{-2\alpha''} (\log p)^{6\alpha''} +\sigma_{\x}^{1/2} .$$
Recall from the discussion before Lemma \ref{lem_mL_p} that $\Bc_p(\mu^{\x},\phi)= \mu^{\x}.$ 
Hence,   the  desired  estimate  follows immediately.
 \hfill $\square$

\medskip

We continue the proof of the main theorem.
We need the following result which is a consequence of  \cite[Lemma 5.3]{BermanBoucksom}.

\begin{lemma}\label{L:BB} 
Consider  a probability measure $\mu$  supported on a compact set $K\subset X$ such that the $L^2$-semi-norm in \eqref{eq:n1.15} is a norm. 
If $S_p$  is  an orthonormal basis of  $H^0(X,L^{p})$ with respect to this norm, then 
the positive measure $\|\det S_p\|^2_{p\phi} \mu^{\otimes N_p}$ is of mass $N_p!$. 
\end{lemma}

\noindent {\bf End of the proof of Theorem \ref{T:abstract}.} Fix a constant $0<\delta<1$ and an orthonormal basis $S_p$ of $H^0(X,L^p)$ with respect to the $L^2$-norm induced by $\mu$ and $\phi$. 
 We first  show that  there is a constant $c>0$  such that for $p\geq 1,$
\begin{equation}\label{e:S_p}
0\leq \log \| \det S_p\|_{L^\infty(K,p\phi)} -  \log \| \det S_p\|_{L^2(\mu,p\phi)}  \leq cN_p  p^{1-\delta}.
\end{equation}
Here, similar to the discussion after Theorem \ref{thm_DMN}, the norm $ \| \det S_p\|_{L^2(\mu,p\phi)}$ is defined using the product probability measure $\mu^{\otimes N_p}$ on $K^{N_p}\subset X^{N_p}$ together with the metric and weight for $(L^p)^{\boxtimes N_p}$, naturally induced by $h_0$ and $\phi$.  

Since $\mu$ is a probability measure, we  have  
 $$\| \det S_p\|_{L^\infty(K,p\phi)}
\geq \| \det S_p\|_{L^2(\mu,p\phi)}.$$
Now, to complete the proof of \eqref{e:S_p}, we only need to show that
\begin{equation}\label{eq_reduction}
\log \| \det S_p\|_{L^\infty(K,p\phi)}
\leq  \log  \| \det S_p\|_{L^2(\mu,p\phi)} 
+O(N_p p^{1-\delta} ).
\end{equation}
 By  \eqref{e:BM}, we get 
$$|s(x)|^2_{p\phi} \leq    A^2 e^{2Ap^{1-\delta}} \| s\|^2_{L^2(\mu,p\phi)}$$
for  every section $s\in H^0(X,L^{p})$, $p\geq 1$, and  
$x\in X.$ If $x_1,\ldots, x_{N_p}$ are points in $X,$ 
then for each $j$
$$x\mapsto \det S_p(x_1,\ldots,x_{j-1}, x, x_{j+1},
\ldots, x_{N_p})$$
is a holomorphic section in $H^0(X,L^p).$ 
A successive application of the last inequality for $j=1,2,\ldots, N_p$ and Fubini's theorem yield
$$\| \det S_p\|^2_{L^\infty(K,p\phi)}\leq  A^{2N_p} 
e^{2AN_pp^{1-\delta}}  \| \det S_p\|^2_{L^2(\mu,p\phi)}.$$
Estimates  \eqref{eq_reduction} and \eqref{e:S_p} follow.

\medskip

Recall that $N_p=O(p^n)$ and by  Stirling's formula $N_p!\approx (N_p/e)^{N_p}\sqrt{2\pi N_p}$. Therefore, 
 Lemma \ref{L:BB} implies
$${1\over p N_p} \log\| \det S_p\|_{L^2(\mu,p\phi)} ={1\over pN_p} \log \sqrt{N_p!} =O(p^{-1}\log p).$$
It follows from  \eqref{e:S_p} that 
\begin{equation}\label{e:S_p_max}
0\leq {1\over pN_p} \log \| \det S_p\|_{L^\infty(K,p\phi)} \leq c_1 p^{-\delta}\quad\text{with some constant}\quad c_1>0.
\end{equation}

Let $\lambda_0>0$ be a constant whose value will be determined later.
For every $p\geq 1$, consider the set
$$ E_p:=\Big\lbrace  \x\in K^{N_p}:\ {1\over pN_p}\log \|\det S_p(\x)\|_{p\phi}\leq - \lambda_0 p^{-\delta}\Big\rbrace.$$
So  for $\x\in K^{N_p}\setminus E_p,$ using  \eqref{e:S_p_max}, we obtain  $\sigma_{\x}\leq (c_1+\lambda_0) p^{-\delta},$ where as above
$$\sigma_{\x}:={1\over pN_p} \log \|\det S_p\|_{L^\infty(K,p\phi)}-{1\over pN_p} \log \|\det S_p(\x)\|_{p\phi}.$$
Hence, 
applying Theorem  \ref{thm_DMN} yields 
\begin{equation*}
\dist_\gamma(\mu^{\x}, \mu_\eq(K,\phi))\leq c p^{-\alpha''\gamma}(\log p)^{3\alpha''\gamma}+c p^{-\gamma\delta/4},
\end{equation*}
for some constant $c>0$. 

To complete the proof of the theorem, it remains to bound the size of $E_p$.
Fix a constant $\lambda$ as in Theorem \ref{T:abstract}. 
Consider two different cases according to the value of $\beta.$

\medskip

\noindent{\bf Case 1.}  Assume that $\beta\geq 2.$ Choose $\lambda_0=\lambda/\beta$. 
We first bound the mass of  $\|\det S_p\|^\beta_{p\phi} \mu^{\otimes N_p}$ from below. Recall that $\mu^{\otimes N_p}$ is a probability measure. Applying
  H\"older's inequality and using   Lemma \ref{L:BB}, we  obtain 
  $$\int \|\det S_p\|^\beta_{p\phi} d\mu^{\otimes N_p}\geq  \Big (\int \|\det S_p\|^2_{p\phi} d\mu^{\otimes N_p} \Big)^{\beta/2}=(N_p!)^{\beta/2}.$$
Consequently,
$\nu^\beta_p\leq \|\det S_p\|^\beta_{p\phi} \mu^{\otimes N_p}. $
Hence, by definition of $E_p$, we get
$$\nu^\beta_p(E_p)\leq \int_{E_p} \|\det S_p (\x)\|^\beta_{p\phi} d\mu^{\otimes N_p}(\x)
\leq \int_{E_p} e^{-\lambda p^{1-\delta} N_p} d\mu^{\otimes N_p}(\x)\\
\leq e^{-\lambda p^{1-\delta} N_p}.$$
This completes the proof for the case $\beta\geq 2$. 

\medskip

\noindent{\bf Case 2. } Assume that $0<\beta\leq 2.$ 
Combining \eqref{e:S_p_max}
and Lemma \ref{L:BB}, we get 
  $$\int_{K^{N_p}}\|\det(S)\|_{p\phi}^\beta d\mu^{\otimes N_p}
\geq  e^{-(2-\beta) c_1p^{1-\delta} N_p}\int_{K^{N_p}}\|\det(S)\|_{p\phi}^2 d\mu^{\otimes N_p}
\geq  e^{-(2-\beta) c_1p^{1-\delta} N_p}.$$  
 Consequently,
$$  \nu^\beta_p\leq e^{(2-\beta) c_1p^{1-\delta} N_p} \|\det S_p\|^\beta_{p\phi} \mu^{\otimes N_p}.
  $$
Hence, we infer 
\begin{eqnarray*}
\nu^\beta_p(E_p)&\leq &  e^{(2-\beta) c_1p^{1-\delta} N_p} \int_{E_p}\|\det S_p (\x)\|^\beta_{p\phi}d \mu^{\otimes N_p}(\x)\\
&\leq&   e^{(2-\beta) c_1p^{1-\delta} N_p} \int_{E_p} e^{-\beta \lambda_0 p^{1-\delta} N_p} d\mu^{\otimes N_p}(\x)\\
&\leq& e^{p^{1-\delta} N_p\big ((2-\beta) c_1 -\beta\lambda_0 \big)}.
\end{eqnarray*}
Choose $\lambda_0 \gg c_1$ and  the result follows. This ends the proof of our main theorem.
 \hfill $\square$

\medskip

As mentioned above, Theorem \ref{thm_DMN} can be applied to other situations. We present now one more application. Consider the same setting as in Theorem \ref{thm_DMN} and a probability measure $\mu$ on $K$.  Recall the following notion, see \cite{BBW}.

\begin{definition} \rm
 Let $0 <r \leq \infty$ and $0<r'\leq\infty$.  
 We say that   $\y\in K^{N_p}$ is  an {\it $(r,r')$-optimal configuration of order $p$}
 if the  following  function in $\x\in K^{N_p}$ 
$$\tau_\x:=\sup_{s\in H^0(X,L^p)\setminus \{0\}} { \|s\|_{L^r(\mu,p\phi)}\over \|s\|_{L^{r'}(\mu^{\x},p\phi)}}$$   
attains its minimum  at $\y.$
\end{definition} 

We have the following elementary property, see also \cite[Proposition 2.10]{BBW}.

\begin{lemma} \label{l:distortion}
If $\y\in K^{N_p}$ is  $(r,r')$-optimal, then $\tau_\y\leq N_p^{1+1/r'}$.
\end{lemma}
\proof
Let $\x=(x_1,\ldots,x_{N_p})$ be a Fekete configuration of order $p$. We only need to check that $\tau_\x\leq N_p^{1+1/r'}$. Choose a basis $S_p=(s_1,\ldots, s_{N_p})$ of $H^0(X,L^p)$ such that $s_i(x_j)=0$ when $i\not=j$ and $\|s_i(x_i)\|_{p\phi}=1$. 
Since $\x$ is a Fekete configuration, we have $\|\det S_p(\cdot)\|_{p\phi}\leq 1$ on $K^{N_p}$. This inequality on $K_i:=\{x_1\}\times\cdots\times \{x_{i-1}\}\times K\times \{x_{i+1}\}\times\cdots\times \{x_{N_p}\}$ implies that $\|s_i(\cdot)\|_{p\phi}\leq 1$ on $K$. 
Finally, if $s$ is a section in $H^0(X,L^p)\setminus \{0\}$, write $s=\lambda_1s_1+\cdots +\lambda_{N_p}s_{N_p}$ and we have
$${\|s\|_{L^r(\mu,p\phi)}\over \|s\|_{L^{r'}(\mu^{\x},p\phi)}}\leq {\sum |\lambda_i| \over (N_p^{-1}\sum |\lambda_i|^{r'})^{1/r'}}\leq {N_p\max |\lambda_i| \over (N_p^{-1}\max |\lambda_i|^{r'})^{1/r'}}=N_p^{1+1/r'}.$$   
The lemma follows.
\endproof

We deduce from Theorem \ref{thm_DMN} the following result, where the simple convergence of $\mu^\y$ when $p\to\infty$ was established in \cite{BBW}.

\begin{corollary} \label{c:Lr-optimal}
In the setting of  Theorem \ref{T:abstract}, consider two  numbers  $0<r, r' \leq \infty$. 
There is a constant $c>0$ such that if $\y$ is an $(r,r')$-optimal configuration of order $p$ for some $p>1$, then   
$$\dist_\gamma(\mu^{\y}, \mu_\eq(K,\phi))\leq c q^{\gamma}.$$
\end{corollary}
\proof
We only have to check that 
$$\sigma_\x\leq c( p^{-1} \log{\tau_\x}+p^{-\delta})\quad\text{for}\quad \x\in K^{N_p}$$
for some constant $c>0$.
Then, Theorem \ref{thm_DMN}, Lemma \ref{l:distortion} and the estimate $N_p=O(p^n)$ imply the result.

We can assume that $\det S_p(\x)\not=0$ because  the case  $\det S_p(\x)=0$ is trivial. So we can choose $S_p=(s_1,\ldots,s_{N_p})$ as in the proof of Lemma \ref{l:distortion}, but here $\x$ is no more a Fekete configuration.
By definition of $\tau_\x$, we have
$$\|s_i\|_{L^r(\mu,p\phi)} \leq \tau_x \|s_i\|_{L^{r'}(\mu^\x,p\phi)}=N_p^{-1/r'} \tau_x\leq \tau_\x.$$
Hence, it follows from Lemma \ref{l:Lr} that 
$$\|s_i\|_{L^\infty(K,p\phi)} \leq A' e^{A'p^{1-\delta}}\tau_x.$$
Therefore, we get
$$\|\det S_p(\cdot)\|_{p\phi}\leq N_p! \big(A' e^{A'p^{1-\delta}} \tau_x\big)^{N_p} \quad \text{on} \quad K^{N_p}.$$
We then deduce the desired estimate using the definition of $\sigma_\x$ and that $\|\det S_p(\x)\|_{p\phi}=1$ by the choice of $S_p$.
\endproof

\end{document}